\documentstyle[preprint,aps]{revtex}

\begin{document}
\title{A brief note on how Einstein's general relativity has influenced the
development of modern differential geometry. }
\author{C. Romero$^{a}$ and F. Dahia$^{b}$}
\address{$^{a}$Departamento de F\'{\i}sica, Universidade Federal da Para%
\'{\i}ba, C.Postal 5008, 58051-970 Jo\~{a}o Pessoa, Pb, Brazil\\
$^{b}$Departamento de F\'{\i}sica, Universidade Federal de Campina Grande,\\
58109-970, Campina Grande, Pb, Brazil\\
E-mail: cromero@fisica.ufpb.br}
\maketitle
\pacs{04.50.+h, 04.20.Cv}

\begin{abstract}
We briefly review a few aspects of the development of differential geometry
which may be considered as being influenced by Einstein%
\'{}%
s general relativity. We focus on how Einsteins's quest for a complete
geometrization of matter and electromagnetism gave rise to an enormous
amount of theoretical work both on physics and mathematics. In connection
with this we also bring to light how recent investigation on theoretical
physics has led to new results on some branches of modern differential
geometry.
\end{abstract}

\section{Introduction}

It is almost impossible to give a fair account of all consequences brought
about by Einstein%
\'{}%
s scientific work on the development of modern human thought. In physics,
Einstein%
\'{}%
s ideas were so revolutionary that pratically no branch of this science has
escaped its influence. Besides physics, \ the newborn concepts of space,
time and space-time, as well as the quantum nature of the microscopic world,
have had great impact also on the field of philosophy. Indeed, as early as
1922, the French philosopher Henri Bergson \cite{Bergson}, who was
originally trained in mathematics, wrote a polemical and critical work on
the notion of time coming from special relativity. After the formulation of
the general theory of relativity and the birth of relativistic cosmology
many other philosophers felt imediately impelled to discuss philosophical
aspects of the new theory. The British philosopher Bertrand Russel is a good
example: his book on the theory of relativity, published in 1925 \cite%
{Russel}\ is a nice and a pedagogical account of Einstein' ideas written
with a logical positivism flavour. The apparent success of non-Euclidean \
geometry to describe our physical world seemed to radically discard the
well-established Kant \cite{Kant} concepts of space and time as {\it a
priori }notions and this issue is still a subject of debate\ among
philosophers of science. Then comes mathematics. The impact Einstein%
\'{}%
s general relativity has had in mathematics is immense. No wonder that it
has not, as far as we know, been fully assessed by historians of science. Of
course we have no intention to embark on such a endeavour here. In this
article, our aim merely consists in pointing out a few particular
mathematical developments which in our view were directly stimulated by
Einstein%
\'{}%
s ideas. In the first section we give a short outline of the general theory
of relativity, some important historical facts and later developments. In
the second section, we a take as a case study the discovery of some
embedding theorems of differential geometry and show how they were
physically motivated in the light of modern theoretical physical research.

\section{The general theory of relativity ( a very brief outline)}

Historically the general theory of relativity (1915) grew out of the special
theory (1905). The mathematical structure of the later in its original
formulation was very simple. However, soon after the appearance of the
special theory, two mathematicians, Hermann Minkowski \cite{Minkowski} and
his contemporary colleague Henri Poincar\'{e} \cite{Poincare} made
significant contributions to its mathematical structure by realizing that
the set of all Lorentz transformations, i.e. those which relate two
different inertial reference frames, constitutes a group and that this group
leaves invariant a certain quadratic form defined in a four-dimensional
space $M^{4}$\ (or Minkowski space-time). This invariant is now referred in
all relativity textbooks as the interval (pseudo-distance) between two
events (points in $M^{4}$) . In view of this discovery it would be more
natural in the context of special relativity theory to\ treat space and time
no more as separate entities, but as mixed together into a new entity, the
space-time. The reaches of this apparently innocuous finding were to be
tremendous. Two comments are in order. First, a new branch of mathematics
was born. Lorentz invariance stimulated the investigation of a new kind of
manifolds endowed with indefinite metrics, now known as semi-Riemannian (or
pseudo-Riemannian) manifolds \cite{ONeill}. Second, from the standpoint of
physics, there was the hint that the new relativistic theory of gravitation
ought to be formulated in a four-dimensional space-time, and that, combined
with the Principle of Equivalence, ultimately led to the geometrization of
the gravitation field, and it is here that lies the astounding beauty of the
general relativity theory. Physics and geometry are identified. and
geometrical curvature mimics the effects of gravitational forces acting on
particles.

General relativity assumes that in the presence of gravitation our
space-time is best represented by a four-dimensional manifold endowed with a
Lorentzian metric. It says nothing about the space-time global topology, so
in this respect it is still a local theory \cite{Luminet}. An elegant set of
partial hyperbolic non-linear equations, found by Einstein and Hilbert \cite%
{Einstein,Hilbert}, is used to determine the metric fields from the
distribution of matter in space-time. Einstein himself did not expect to
solve his field equations exactly and the first solution was obtained by
Karl Schwarzschild in 1916 \cite{Schwarzschild}. Schwarzschild%
\'{}%
s solution describes the geometry of the space-time outside a spherically
symmetric matter distribution and contained two puzzling features: the
existence of an event horizon and a space-time singularity. Both these
aspects of Schwarzschild's solution, which is the prototype of a noncharged
static blackhole, were to generate a great deal of mathematical work in the
following years.

Very soon general relativity theory was applied to cosmology. In 1917
Einstein wrote a paper in which he modifies the field equations to tackle
the problem of finding the geometry of the Universe \cite{Einstein2} . His
cosmological model described a homogeneous, isotropic and static universe
whose spatial geometry may be viewed as the geometry of a hypersphere
embedded in an Euclidean four-dimensional space. This was a nice example of
a finite universe with no boundaries. However, Einstein's universe did not
account for the recession motion of galaxies, observed in 1929 by the
American astronomer Edwin Hubble. This discovery of this effect, which was
interpreted by the Belgian physicist Georges Lema\^{\i}tre \cite{Lemaitre}
as an evidence of the expansion of the Universe, would drastically change
our view of the Cosmos. Indeed, the only plausible explanation of the fact
that galaxies are moving away from us is that the Universe is expanding.
Curiously enough, an expanding solution of Einstein's original field
equations had already been obtained by a little known Russian scientist,
Alexander Friedmann, in 1922. \ Friedmann's time-dependent solution
introduced a revolutionary ingredient in our view of the Universe: the idea
that the Cosmos started out with a big bang. In mathematical terms it means
that the geometry of Friedmann's model, like Schwarzschild's solution,
contains a singularity (space-time is geodesically incomplete ). Careful
investigation of the nature and mathematical structure of singularities
found in solutions of Einstein's field equations ultimately led Roger
Penrose \cite{Penrose} and Stephen Hawking \cite{Hawking}, in the sixties,
to discover the famous singularities theorems which have strongly boosted
the study of global aspects of general relativity \cite{Geroch} where
methods of differential topology have been extensively employed to
investigate the problems \cite{Penrose1}.

\section{General relativity and differential geometry}

One of the most cherished projects of contemporary theoretical physicists is
to find a theory capable of unifying the fundamental forces of nature, a
theory of everything, as it has been called. Unification, in fact, has been
a feature of all great theories of physics. In a certain sense Newton,
Maxwell and Einstein, they all succeeded in performing some sort of
unification. Twentieth century physics has recurrently pursued this theme.
Now broadly speaking one can mention two different paths followed by
theoreticians to arrive at unified field theory.\ First there are the early
attempts of Einstein, Weyl, Cartan, Eddington, Schr\"{o}dinger and many
others, whose task consisted of unifying gravity and electromagnetism \cite%
{Goenner}. The methodological approach of this group consisted basically in
resorting to different kind of non-Riemannian geometries capable of
accomodating new geometrical structures with a sufficient number of degrees
of freedom to describe the electromagnetic field. In this way different
types of geometry have been "created", such as affine geometry (asymetric
connection), Weyl's geometry (where the notion of parallel transport differs
from Levi-Civita's notion), etc. It is not easy to track further
developments of these geometries motivated by general relativity. Already in
1921 the Dutch mathematician Schouten wrote: \textquotedblleft Motivated by
relativity theory, differential geometry received a totally novel, simple
and satisfying foundation" (quoted in \cite{Goenner}). However, the snag
with all these attempts was that they completely ignored quantum mechanics
and dealt with unification only in a classical level. Of course, an approach
to unification today would necessarily take into account quantum field
theory. Now the second approach to unification comes into play. It has to do
with the rather old idea that our space-time may have more than four
dimensions.

The story starts with the work of the Finnish physicist Gunnar Nordstr\"{o}m 
\cite{Nordström}, in 1914. Nordstr\"{o}m realised that by postulating the
existence of a fifth dimension he was able (in the context of his scalar
theory of gravitation) to unify gravity and electromagnetism by embedding
space-time into a five-dimensional space. Although the idea was quite
original and interesting it seems the paper did not attract much attention
due to the fact that his gravitation theory was not accepted at the time.
Then, soon after the completion of general relativity, Th\'{e}odor Kaluza,
and later, Oscar Klein, launched again the same idea, now entirely based on
Einstein's theory of gravity. In a very creative manner the Kaluza-Klein
theory \ starts from five-dimensional vacuum Einstein's equations and show
that, under certain assumptions, they reduce to a four-dimensional system of
coupled Einstein-Maxwell equations. The paper was seminal and gave rise to
several different theoretical developments exploring the idea of achieving
unification from extra dimensionality of space. Indeed, through the old and
modern versions of Kaluza-Klein theory \cite{kaluza,kaluza1,kaluza2},
supergravity \cite{deser}, superstrings \cite{superstrings}, and to the more
recent braneworld scenario\cite{rs1,rs2}, induced-matter \cite{overduin,book}
and M-theory \cite{duff}, there is a strong belief among some physicists
that unification might be finally achieved if one accepts that space-time
has more than four dimensions.

Amidst all these higher-dimensional theories, one of them, the
induced-matter theory (also referred to as space-time-matter theory \cite%
{overduin,book}) has called our attention for it recalls Einstein%
\'{}%
s belief that matter and radiation (not only the gravitation field) should
be viewed as manifestations of pure geometry\cite{Einstein3}. Kaluza-Klein
theory was a first step in this direction. But it was Paul Wesson \cite{book}%
, from the University of Waterloo, who pursued the matter further. Wesson
and collaborators realized that by embedding the ordinary space-time into a
five-dimensional vacuum space, it was possible to describe the macroscopic
properties of matter in geometrical terms. In a series of interesting papers
Wesson and his group showed how to produce standard cosmological models from
five-dimensional vacuum space. It looked like any energy-momentum tensor
could be generated by an embedding mechanism. At the time these facts were
discovered, there was no guarantee that {\it any }energy-momentum could be
obtained in this way. Putting it in mathematical terms, Wesson's programm
would not work always unless one could prove that {\it any }solution of
Einstein's field equations could be isometrically embedded in
five-dimensional Ricci-flat space \cite{Romero}. As it happens, that was
exactly the content of a beautilful and powerful theorem of differential
geometry now known as the Campbell-Magaard theorem \cite{Campbell}. Although
very little known, the theorem was articulated by the English mathematician
John Campbell in 1926 and was given a complete proof only in 1963 by Lorenz
Magaard \cite{magaard}. (At this point may we digress a little bit.
Campbell, who died in 1924 \cite{Campbell1}, \ was interested in geometrical
aspects of Einstein's relativity and his works \cite{Campbell2} were
published a few years before the classical Janet-Cartan \cite{janet,cartan}
theorem on embeddings was established \footnote{%
Janet-Cartan theorem originated from a conjecture by Schl\"{a}ff, in 1873,
and states that if the embedding space is flat, then the maximum number of
extra dimensions needed to analytically embed a Riemannian manifold is $d$ ,
with $0\leq d\leq n(n-1)/2$. The novelty brought by Campbell-Magaard theorem
is that the number of extra dimensions falls drastically to $d=1$ when the
embedding manifold is allowed to be Ricci-flat (instead of Riemann-flat).}.
Manifolds called {\it Einstein spaces} had begun to attract the interest of
mathematicians soon after the discovery of Schwarzschild space-time and
de-Sitter cosmological models). Now compared to the Janet-Cartan theorem the
nice thing about the Campbell-Magaard's result is that the codimension of
the embedding space is drastically reduced: one needs only one
extra-dimension, and that perfectly fits the requirements of the
induced-matter theory. Finally, let us note both theorems refer to local and
analytical embeddings (the global version of Janet-Cartan theorem was worked
out by John Nash \cite{Nash}, in 1956, and adapted for semi-Riemannian
geometry by R. Greene \cite{Greene}, in 1970,\ while a discussion of global
aspects of Cambell-Magaard has recently appeared in the literature \cite{Kat}%
).

\section{Higher-dimensional space-times and the search for new theorems}

Apart from induced-matter theory, there appeared at the turn of the XX
century some other physical models of the Universe, which soon attracted the
attention of theoreticians. These models have put forward the idea that
ordinary space-time may be viewed as a four-dimensional hypersurface
embedded not in a Ricci-flat space, but in a five-dimensional Einstein space
(referred to as {\it the bulk}) \cite{Randall}. Spurred by this proposal new
research on the geometrical structure of the proposed models started. It was
conjectured \cite{anderson} and later proved that the Campbell-Magaard could
be immediately generalized for embedding Einstein spaces \cite{dahia1} This
was the first extension of the Campbell-Magaard theorem and other extensions
were to come. More general local isometric embeddings were next
investigated, and it was proved that any $n$-dimensional semi-Riemannian
analytic manifold can be locally embedded in \ $(n+1)$-dimensional analytic
manifold with a non-degenerate Ricci-tensor, which is equal, up to a local
analytical diffeomorphism, to the Ricci-tensor \ of an arbitrary specified
space \cite{dahia3}. Further motivation in this direction came from studying
embeddings in the context of non-linear sigma models, a theory proposed by
J. Schwinger in the fifties to describe strongly interacting massive
particles \cite{Schwinger}. It was then showed that any $n$-dimensional
Lorentzian manifold $(n\geq 3)$ can be harmonically embedded in a $(n+1)$%
-dimensional semi-Riemannian Ricci-flat manifold \cite{dahia4}. As a final
remark on the Campbell-Maggard theorem and its application to physics, let
us note that its proof is based on the Cauchy-Kowalevskaya theorem.
Therefore, some properties of relevance to physics, such as the stability of
the embedding, cannot be guaranteed to hold \cite{Anderson}. Nevertheless,
the problem of embedding space-time into five-dimensional spaces can be
considered in \ the context of the Cauchy problem \ in general relativity 
\cite{Yvonne}. Specifically, it has recently been shown that the embedded
space-time may arise as a result of physical evolution of proper initial
data. This new perspective has some advantages in comparison with the
original Campbell-Magaard formulation because it allows us, by exploring the
hyperbolic character of Einstein field equations, to show that the embedding
has stability and domain of dependence (causality) properties \cite{dahia5}.

\section{Conclusion}

We would like to conclude by pointing out that all these developments
essentially grew out of one great theory: General relativity. Underlying
this connection between physics and geometry there is the basic idea that a
theory of the gravitational field must be a {\it metric} theory. Now there
is a vast number of metric theories. Their motivation is twofold:
quantization of gravity and its unification with the other physical fields.
Some of these theories postulate the existence of extra dimensions of the
Universe and these multidimensional theories of space-time have employed a
complex and sophisticated mathematical language, imported from modern
differential geometry and topology. That strange belief on "the unreasonable
effectiveness of mathematics in the natural sciences", as put by the
physicist E. Wigner \cite{Wigner} , seems to be still alive among
contemporary physicists. However, this is a mutual process of interaction
between the two sciences. In this paper we have tried to explore the other
side of this relationship, and how physical research can be beneficial to
the development of mathematics itself, in particular the important role
Einstein's general relativity has played in promoting progress of some
branches of modern differential geometry.

\bigskip

\bigskip

\section{Acknowledgement}

The authors thank CNPq for financial support.

\end{document}